\documentclass[11pt]{article}
\usepackage{graphicx}
\usepackage{amsmath,amsfonts,amssymb}
\usepackage{url}
\usepackage{amsmath}
\usepackage{graphicx}
\usepackage{epsfig}
\usepackage{epstopdf}
\usepackage{color}
\def\tto{\;{\lower 1pt \hbox{$\rightarrow$}}\kern -10pt
\hbox{\raise 2pt \hbox{$\rightarrow$}}\;}
\def\Hat{\widehat}

\def\Bar{\overline}
\def\ra{\rangle}
\def\la{\langle}
\def\ve{\varepsilon}
\def\epsilon{\varepsilon}
\def\B{\Bbb B}
\def\h{\hfill\Box}
\def\R{\Bbb R}

\def\N{I\!\!N}
\def\ox{\bar{x}}

\def\core{\mbox{\rm core}\,}

\def\dom{\mbox{\rm dom}\,}

\def\i{\mbox{\rm int}\,}

\def\i{\mbox{\rm int}\,}

\def\h{\hfill\square}
\def\dn{\downarrow}
\def\O{\Omega}
\def\ph{\varphi}
\def\emp{\emptyset}
\def\st{\stackrel}
\def\oR{\Bar{\R}}

\def\gg{\gamma}
\def\bb{\beta}
\def\dd{\delta}
\def\al{\alpha}
\def\Th{\Theta}
\def\N{I\!\!N}

\def\Lm{\Lambda}
\newcounter{lk}

\begin{document}
\begin{center}
{\bf EXTREMALITY OF CONVEX SETS WITH SOME APPLICATIONS}\\[3ex]
BORIS S. MORDUKHOVICH\footnote{Corresponding author. Department of Mathematics, Wayne State University, Detroit, MI 48202, USA (boris@math.wayne.edu) and Peoples' Friendship University of Russia, Moscow 117198, Russia. Email: boris@math.wayne.edu, phone: (734)369-3675, fax: (313)577-7596. Research of this author was partly supported by the National Science Foundation under grants DMS-1007132 and DMS-1512846 and by the Air Force Office of Scientific Research under grant \#15RT0462.} and NGUYEN MAU NAM\footnote{Fariborz Maseeh Department of Mathematics and Statistics, Portland State University, Portland, OR 97207, USA. Email: mau.nam.nguyen@pdx.edu. Research of this author was partly supported by the National Science Foundation under grant \#1411817.}.\\[3ex]
{\bf Dedicated to the memory of Jonathan Michael Borwein}
\end{center}
\small{\bf Abstract:} In this paper we introduce an enhanced notion of extremal systems for sets in locally convex topological vector spaces and obtain efficient conditions for set extremality in the convex case. Then we apply this machinery to deriving new calculus results on intersection rules for normal cones to convex sets and on infimal convolutions of support functions.\\[1ex]

\newtheorem{Theorem}{Theorem}[section]
\newtheorem{Proposition}[Theorem]{Proposition}
\newtheorem{Remark}[Theorem]{Remark}
\newtheorem{Lemma}[Theorem]{Lemma}
\newtheorem{Corollary}[Theorem]{Corollary}
\newtheorem{Definition}[Theorem]{Definition}
\newtheorem{Example}[Theorem]{Example}
\renewcommand{\theequation}{\thesection.\arabic{equation}}
\normalsize\vspace*{-0.2in}

\newpage
\section{Introduction}\vspace*{-0.1in}

{\em Convex analysis} has been well recognized as an important area of mathematics with numerous applications to optimization, control, economics, and many other disciplines. We refer the reader to the fundamental monographs \cite{bc,bl,HU,r,z} and the bibliographies therein for various aspects of convex analysis and its applications. Jon Borwein, who unexpectedly passed away on August 2, 2016, made pivotal contributions to these and related fields of Applied Mathematics, among other areas of his fantastic creative activity.

Methods and constructions of convex analysis play also a decisive role in the study of nonconvex functions and sets by using certain convexification procedures. In particular, calculus and applications of Clarke's generalized gradients for nonconvex functions \cite{c} is based on appropriate convexifications and employing techniques and results of convex analysis.

Besides this, other ideas have been developed in the study and applications of nonconvex functions, sets, and set-valued mappings in the framework of {\em variational analysis}, which employs variational/optimization principles married to perturbation and approximation techniques; see the books \cite{BZ,m-book1,RockWets-VA} for extended expositions in finite and infinite dimensions. Powerful tools, results, and applications of variational analysis have been obtained by using the {\em dual-space geometric approach} \cite{m-book1} based on the {\em extremal principle} (a geometric variational principle) for systems of sets. This approach produces first a {\em full calculus} of generalized normals to nonconvex sets and then applies it to establish comprehensive calculus rules for related subgradients of extended-real-valued functions and coderivatives of set-valued mappings. Needless to say that well-developed calculus of generalized differentiation is an unavoidable requirement and the key for various applications.

Addressing generally nonconvex objects, results of variational analysis contain corresponding convex facts as their particular cases. However, basic variational techniques involving limiting procedures do not fully capture advantages from the presence of convexity. Indeed, the major calculus results of \cite{m-book1} hold in {\em Asplund} spaces (i.e., such Banach spaces where every separable subspace has a separable dual) and the {\em closedness} of sets (epigraphs for extended-real-valued function, graphs for set-valued mappings) is a standing assumption.

The major goal of this paper is to investigate a counterpart of the variational geometric approach to the study of convex sets in locally convex topological vector (LCTV) spaces without any completeness and closedness assumptions. Based on an enhanced notion of {\em set extremality}, which is a global version of the corresponding local concept largely developed and applied in \cite{m-book1} while occurring to be particularly useful in the convex setting mainly exploited here, this approach allows us to obtain the basic intersection rule for normals to convex sets under a new qualification condition. The same approach also allows us to derive new calculus results for support functions of convex set intersections in general LCTV spaces. Note that these results can be used to obtain major calculus rules of generalized differentiation and Fenchel conjugates for extended-real-valued convex functions; cf.\ our previous publications \cite{bmn,bmn1} for some versions in finite dimensions.

The rest of the paper is organized as follows. In Section~2 we introduce the aforementioned version of set extremality, establish its relationships with the separation property for convex sets, and derive various extremality conditions. The obtained results are applied in Section~3 to get the normal cone representation for convex set intersections under a new qualification condition. In Section~4 this approach is employed to represent the support function of set intersections via the infimal convolution of supports to intersection components.

For simplicity of presentation we suppose, unless otherwise stated, that all the spaces under consideration are {\em normed linear} spaces. The reader can check that the results obtained below in this setting hold true in the LCTV space generality.

The notation used throughout the paper is standard in the areas of functional, convex, and variational analysis; cf.\ \cite{m-book1,r,RockWets-VA,z}. Recall that the closed ball centered at $\ox$ with radius $r>0$ is denoted by $\B(\ox;r)$ while the closed unit ball of the space $X$ in question and its topological dual $X^*$ are denoted by $\B$ and $\B^*$, respectively, if no confusion arises. Given a convex set $\O\subset X$, we write $\R^+(\Omega):=\{tv\in X|\;t\in\R_+,\;v\in\Omega\}$, where $\R_+$ signifies the collection of positive numbers, and use the symbol $\overline\O$ for the topological closure of $\O$. Finally, remind the notation for the (algebraic) {\em core} of a set:
\begin{equation}\label{core-def}
\core\Omega:=\big\{x\in\Omega\big|\;\forall\,v\in X\;\exists\,\gg>0\;\mbox{\rm such that }\;x+tv\in\Omega\;\mbox{\rm whenever }\;|t|<\gg\big\}.
\end{equation}
In what follow we deal with {\em extended-real-valued} functions $f\colon X\to\oR:=(-\infty,\infty]$ and assume that are {\em proper}, i.e., $\dom f:=\{x\in X|\;f(x)<\infty\}\ne\emp$.\vspace*{-0.2in}

\section{Extremal Systems of Sets}
\setcounter{equation}{0}\vspace*{-0.1in}

We start this section with the definition of extremality for set systems, which is inspired by the notion of local set extremality in variational analysis (see \cite[Definition~2.1]{m-book1}) while having some special features that are beneficial for convex sets. In particular, we do not require that the sets have a common point.\vspace*{-0.1in}

\begin{Definition}{\bf(set extremality).}\label{ext-sys} We say that two nonempty sets $\O_1,\O_2\subset X$ form an {\sc extremal system} if for any $\epsilon>0$ there exists $a\in X$ such that
\begin{equation}\label{setex}
\|a\|\le\epsilon\;\;\mbox{\rm and }\;(\Omega_1+a)\cap\Omega_2=\emp.
\end{equation}
\end{Definition}\vspace*{-0.05in}

Observe similarly to \cite{m-book1} that the notion of set extremality introduced in Definition~\ref{ext-sys} covers (global) optimal solutions to problems of constrained optimization with scalar, vector, and set-valued objectives, various equilibrium concepts arising in operations research, mechanics, and economic modeling, etc. Furthermore, the set extremality naturally arises in deriving calculus rules of generalized differentiation in variational analysis. In particular, we are going to demonstrate this below in our device of the normal cone intersection rule and the support function representation for convex set intersections presented in the paper.

Given a convex set $\Omega\subset X$ with $\ox\in\Omega$, the {\em normal cone} to $\Omega$ at $\ox$ is
\begin{equation}\label{nor}
N(\ox;\Omega):=\big\{x^*\in X^*\big|\;\la x^*,x-\ox\ra\le 0\;\;\mbox{\rm for all }\;x\in\Omega\big\}.
\end{equation}

The following underlying result establishes a useful characterization of set extremality and shows that, in the case of convex sets, extremality is closely related to while being different from the conventional convex separation:
\begin{equation}\label{sep}
\sup_{x\in\O_1}\la x^*,x\ra\le\inf_{x\in\O_2}\la x^*,x\ra\;\;\mbox{\rm for some }\;x^*\ne 0.
\end{equation}
Note that if $\O_1, \O_2$ are convex sets  such that $\ox\in\O_1\cap\O_2$, then \eqref{sep} is equivalent to
\begin{eqnarray}\label{ep}
N(\ox;\O_1)\cap\big(-N(\ox;\O_2)\big)\ne\{0\}.
\end{eqnarray}\vspace*{-0.35in}
\begin{Theorem}{\bf(set extremality and separation).}\label{extremal principle} Let $\O_1,\O_2\subset X$ be nonempty sets. Then the following assertions are fulfilled:

{\bf(i)} The sets $\O_1$ and $\O_2$ form an extremal system if and only if $0\notin{\rm int}(\Omega_1-\Omega_2)$. Furthermore, the extremality of $\O_1,\O_2$ implies that $({\rm int}\,\O_1)\cap\O_2=\emp$ and likewise $({\rm int}\,\O_2)\cap\O_1=\emp$.

{\bf(ii)} If $\O_1,\O_2$ are convex and form an extremal system and if ${\rm int}(\O_1-\O_2)\ne\emp$, then the separation property \eqref{sep} holds.

{\bf (iii)} The separation property \eqref{sep} always implies the set extremality \eqref{setex}, without imposing either the convexity of $\O_1,\O_2$ or the condition ${\rm int}(\O_1-\O_2)\ne\emp$ as in {\rm(ii)}.
\end{Theorem}\vspace*{-0.1in}
{\bf Proof.} To verify the extremality characterization in (i), suppose first that the sets $\O_1,\O_2$ form an extremal system while the condition $0\notin{\rm int}(\Omega_1-\Omega_2)$ fails. Then there is $r>0$ such that $\B(0;r)\subset\O_1-\O_2$. Put $\ve:=r$ and observe that $-a\in\O_1-\O_2$ for any $a\in X$ with $\|a\|\le\ve$, which gives us $(\O_1+a)\cap\O_2\ne\emp$ and thus contradicts \eqref{setex}. To justify the converse implication in (i), suppose that $0\notin{\rm int}(\Omega_1-\Omega_2)$. Then for any $\ve>0$ we get
$$
\B(0;\ve)\cap\big(X\setminus(\O_1-\O_2)\big)\ne\emp,
$$
which tells us that there is $a\in X$ such that $\|a\|<\ve$ and $-a\in\O_1-\O_2$, i.e., \eqref{setex} holds. It remains to show in (i) that the extremality of $\O_1,\O_2$ yields $({\rm int}\,\O_1)\cap\O_2=\emp$. Assuming the contrary, take $x\in{\rm int}\,\O_1$ with $x\in\O_2$ and find $\ve>0$ such that $x-a\in\O_1$ for any $a\in X$ with $\|a\|<\ve$. This clearly contradicts \eqref{setex} and thus completes the proof of (i).

Next we verify (ii). Consider the two convex sets $\Lm_1:=\O_1-\O_2$ and $\Lm_2:=\{0\}$ in $X$. By the extremality of $\O_1,\O_2$ we have due to (i) that $({\rm int}\,\Lm_1)\cap\Lm_2=\emp$, where $\i\Lm_1\ne\emp$ by the assumption in (ii). The classical separation theorem applied to $\Lm_1,\Lm_2$ tells us that $\sup_{x\in\O_1-\O_2}\la x^*,x\ra\le 0$, which is clearly equivalent to \eqref{sep}. Thus assertion (ii) is justified.

To prove the final assertion (iii), take $x^*\ne 0$ from \eqref{sep} and find $c\in X$ such that $\la x^*,c\ra>0$. For any $\ve>0$ we can select $a:=-c/k$ satisfying $\|a\|<\ve$ when $k\in\N$ is sufficiently large. Let as show that \eqref{setex} holds with this vector $a$. If it is not the case, then there exists $\Hat x\in\O_2$ such that $\Hat x-a\in\O_1$. By the separation property \eqref{sep} we have
$$
\la x^*,\Hat x-a\ra\le\sup_{x\in\O_1}\la x^*,x\ra\le\inf_{x\in\O_2}\la x^*,x\ra\le\la x^*,\Hat x\ra,
$$
which gives us by the above construction of $a\in X$ that
$$
\la x^*,\Hat x\ra-\la x^*,a\ra=\la x^*,\Hat x\ra+k\la x^*,c\ra\le\la x^*,\Hat x\ra,
$$
and therefore $\la x^*,c\ra\le 0$. It contradicts the choice of $c\in X$ and hence justifies assertion (iii) while completing in this way the proof of the theorem. $\h$\vspace*{-0.1in}

\begin{Corollary}{\bf (sufficient conditions for extremality of convex sets).}\label{int-ext} Let $\O_1,\O_2$ be nonempty convex sets of $X$ satisfying the conditions $\mbox{\rm int}\,\Omega_1\ne\emp$ and $(\mbox{\rm int}\,\Omega_1)\cap\Omega_2=\emp$. Then the sets $\Omega_1$ and $\Omega_2$ form an extremal system. Furthermore, we have $\mbox{\rm int}(\Omega_1-\Omega_2)\ne\emp$.
\end{Corollary}\vspace*{-0.1in}
{\bf Proof.} It is well known that the assumptions imposed in the corollary ensure the separation property for convex sets. Thus the set extremality of $\O_1,\O_2$ follows from Theorem~\ref{extremal principle}(iii). To verify the last assertion of the corollary, take any $\ox\in{\rm int}\,\O_1$ and find $r>0$ such that ${\rm int}\,\B(\ox;r)\subset\Omega_1$. Then for any fixed point $x\in\O_2$ we have
$$
V:={\rm int}\,\B(\ox;r)-x\subset\Omega_1-\Omega_2,
$$
and thus $\mbox{\rm int}(\Omega_1-\Omega_2)\ne\emp$ because $V$ is a nonempty open subset of $X$. $\h$\vspace*{-0.1in}

\begin{Remark}{\bf (on the extremal principle).}\label{ext-prin} {\rm Condition \eqref{ep} is known to hold, under the name of the (exact) {\em extremal principle}, for locally extremal points of nonconvex sets. In \cite[Theorem~2.22]{m-book1} it is derived for closed subsets of Asplund spaces with the replacement of \eqref{nor} by the basic/limiting normal cone of Mordukhovich, which reduces to \eqref{nor} for convex sets. Besides the Asplund space requirement, the aforementioned result of \cite{m-book1} imposes the {\em sequential normal compactness} (SNC) assumption on one of the sets $\O_1,\O_2$. This property is satisfied for convex sets under the interiority assumption of Corollary~\ref{int-ext}; see \cite[Proposition~1.25]{m-book1}. Furthermore, in the case of closed convex sets in Banach spaces the SNC property offers significant advantages for the validity of \eqref{ep} in comparison with the interiority condition due to the SNC characterization from \cite[Theorem~1.21]{m-book1}: a closed convex set $\O$ with nonempty relative interior (i.e., the interior of it with respect to its span) is SNC at every $\ox\in\O$ if and only if the closure of the span of $\O$ is of finite codimension. A similar characterization has been obtained in \cite[Theorem~2.5]{blm} for the more restrictive Borwein-Str\'ojwas' {\em compactly epi-Lipschitzian} (CEL) property \cite{BS} of closed convex sets in normed spaces. Note that the CEL and SNC properties may not agree even for closed convex cones in nonseparable Asplund spaces; see \cite{fm} for comprehensive results and examples.}
\end{Remark}\vspace*{-0.05in}

As established in Theorem~\ref{extremal principle}(ii), the set extremality in \eqref{setex} implies the separation property \eqref{sep} and its equivalent form \eqref{ep} whenever $\ox\in\O_1\cap\O_2$ under the {\em nonempty difference interior} ${\rm int}(\O_1-\O_2)\ne\emp$ for arbitrary convex sets $\O_1,\O_2$ in LCTV spaces. Could we relax this assumption? The next theorem shows that it can be done, for {\em closed} convex subsets of {\em Banach} spaces, in both {\em approximate} and {\em exact} forms of the {\em convex extremal principle}. Furthermore, the results obtained therein justify that both of these forms are {\em characterizations} of the convex set extremality under the SNC property of one of the sets involved without imposing any interiority assumption on them or their difference.

To proceed, recall first the definition of the SNC property used below for convex sets; compare it with a nonconvex counterpart from \cite[Definition~1.20]{m-book1}. A subset $\O\subset X$ of a Banach space is {\em SNC} at $\ox\in\O$ if for any sequence $\{(x_k,x^*_k)\}_{k\in\N}\subset X\times X^*$ we have
\begin{equation}\label{snc}
\big[x^*_k\in N(x_k;\O),\;x_k\in\O,\;x_k\to\ox,\;x^*_k\st{w^*}{\to}0\big]\Longrightarrow\|x^*_k\|\to 0\;\;\mbox{\rm as }\;k\to\infty,
\end{equation}
where the normal cone is taken from \eqref{nor}, and where the symbol $\st{w^*}{\to}$ signifies the {\em sequential} convergence in the weak$^*$ topology of $X^*$. We have already mentioned in Remark~\ref{ext-prin} the explicit description of the SNC property for closed convex sets with nonempty relative interiors in Banach spaces given in \cite[Theorem~1.21]{m-book1}. Assertion (ii) of the next theorem employs SNC \eqref{snc} for furnishing the limiting procedure in general Banach spaces.\vspace*{-0.1in}

\begin{Theorem}{\bf(approximate and exact versions of the convex extremal principle in Banach spaces).}\label{convex-ep} Let $\O_1$ and $\O_2$ be closed convex subsets of a Banach space $X$, and let $\ox$ be any common point of $\O_1,\O_2$. Consider the following assertions:

{\bf (i)} The sets $\O_i$, $i=1,2$, form an extremal system in $X$.

{\bf (ii)} For each $\ve>0$ we have:
\begin{eqnarray}\label{ep1}
\exists\,x_i\in\B(\ox;\ve)\cap\O_i,\;\exists\,\;x^*_i\in N(x_{i\ve};\O_i)+\ve\B^*\;\;\mbox{\rm with }\;x^*_1+x^*_2=0,\;\|x^*_1\|=\|x^*_2\|=1.
\end{eqnarray}

{\bf (iii)} The equivalent properties \eqref{sep} and \eqref{ep} are satisfied.

Then we always have the implication {\rm(i)}$\Longrightarrow${\rm(ii)}. Furthermore, all the properties in {\rm(i)}--{\rm(iii)} are equivalent if in addition
either $\O_1$ or $\O_2$ is $SNC$ at $\ox$.
\end{Theorem}\vspace*{-0.1in}
{\bf Proof.} Let us begin with verifying (i)$\Longrightarrow$ (ii). It follows from the extremality condition that for any $\epsilon>0$ there exists $a\in X$ such that
\begin{equation*}
\|a\|\le\ve^2\;\;\mbox{\rm and }\;(\Omega_1+a)\cap\Omega_2=\emp.
\end{equation*}
Define the convex, lower semicontinuous, and bounded from below function $f\colon X^2\to\oR$ by
\begin{equation}\label{ext1}
f(x_1,x_2):=\|x_1-x_2+a\|+\dd\big((x_1,x_2);\O_1\times\O_2\big),\quad(x_1,x_2)\in X^2,
\end{equation}
via the indicator function of the closed set $\O_1\times\O_2$. It follows from \eqref{setex} that $f(x_1,x_2)>0$ on $X\times X$ and $f(\ox,\ox)=\|a\|\le\ve^2$ for any $\ox\in\O_1\times\O_2$. Applying to \eqref{ext1} the Ekeland variational principle (see, e.g., \cite[Theorem~2.26(i)]{m-book1}), we find a pair $(x_{1\ve},x_{2\ve})\in\O_1\times\O_2$ satisfying $\|x_{1\ve}-\ox\|\le\ve$, $\|x_{2\ve}-\ox\|\le\ve$, and
\begin{equation*}
f(x_{1\ve},x_{2\ve})\le f(x_1,x_2)+\ve\big(\|x_1-x_{1\ve}\|+\|x_2-x_{2\ve}\|\big)\;\;\mbox{\rm for all }\;(x_1,x_2)\in X^2.
\end{equation*}
The latter means that the function $\ph(x_1,x_2):=f(x_1,x_2)+\ve\big(\|x_1-x_{1\ve}\|+\|x_2-x_{2\ve}\|)$ attains its minimum on $X^2$ at $(x_{1\ve},x_{2\ve})$ with $\|x_{1\ve}-x_{2\ve}-a\|\ne 0$. Thus the generalized Fermat rule tells us that $0\in\partial\ph(x_{1\ve},x_{2\ve})$. Taking into account the summation structure of $f$ in \eqref{ext1}, we apply to its subdifferential the classical Moreau-Rockafellar theorem that allows us to find---by standard subdifferentiation of the norm and indicator functions---such dual elements $x^*_{i\ve}\in N(x_{i\ve};\O_i)+\ve\B^*$ for $i=1,2$ that all the conditions in \eqref{ep1} are satisfied. This justifies assertion (ii) of the theorem.

We verify next the validity of (ii)$\Longrightarrow$(iii) by furnishing the passage to the limit in \eqref{ep1} as $\ve\dn 0$ with the help of the SNC property of, say, the set $\O_1$ at $\ox$. Take a sequence $\ve_k\dn 0$ as $k\to\infty$ and find by \eqref{ep1} the corresponding septuples $(x_{1k},x_{2k},x^*_k,x^*_{1k},x^*_{2k},e^*_{1k},e^*_{2k})$ so that $x_{1k}\to\ox$, $x_{2k}\to\ox$ as $k\to\infty$, and
\begin{equation}\label{ep2}
x^*_k=x^*_{1k}+\ve_ke^*_{1k},\;x^*_k=-x^*_{2k}+\ve_k e^*_{2k},\;\|x^*_k\|=1,\;x^*_{ik}\in N(x_{ik};\O_1),\;e^*_{ik}\in\B^*
\end{equation}
for all $k\in\N$ and $i=1,2$. The classical Banach–-Alaoglu theorem of functional analysis tells us that for any Banach space $X$ the sequence of triples $(x^*_k,e^*_{1k},e^*_{2k})$ contains a {\em subnet} converging to some $(x^*,e^*_1,e^*_2)\in\B^*\times\B^*\times\B^*$ in the weak$^*$ topology of $X^*$. It follows from \eqref{ep2} and definition \eqref{nor} of the normal cone to convex sets that the corresponding subnets of $\{x^*_{1k},x^*_{2k})\}$ converge in the latter topology to some pair $(x^*_1,x^*_2)\in X^*\times X^*$ satisfying $x^*_1=-x^*_2=x^*$ and $x^*_i\in N(\ox;\O_i)$ for $i=1,2$.

To justify (iii), it remains to show that we can always find $x^*\ne 0$ in this way provided that $\O_1$ is SNC at $\ox$. Assuming the contrary, let us first check that $\{x^*_{1k}\}$ converges to zero in the weak$^*$ topology. If it is not the case, there is $z\in X$ such that the numerical sequence $\{\la x^*_{1k},z\ra\}$ does not converge to zero. Fix $w\in \O_1$ and for each $k\in\N$ consider the set
\begin{equation}\label{Vk}
V_k:=\big\{z^*\in X^*\big|\;|\la z^*,w-\ox\ra-\la x^*_1,w-\ox\ra|<1/k,\;|\la z^*,z\ra-\la x^*_1,z\ra|<1/k\big\},
\end{equation}
which is a neighborhood of $x^*_1$ in the weak$^*$ topology of $X^*$. By extracting numerical subsequences in \eqref{Vk}, suppose without loss of generality that
\begin{equation*}
\la x^*_{1k},w-\ox\ra\to\la x^*_1,w-\ox\ra\;\;\mbox{\rm and }\;\la x^*_{1k},z\ra\to\la x^*_1,z\ra\;\;\mbox{\rm as }\;k\to\infty.
\end{equation*}
Remembering that $x^*_{1k}\in N(x_{1k};\O_1)$ by \eqref{ep2} gives us the estimate
\begin{equation}\label{xk}
\la x^*_{1k},w-\ox\ra=\la x^*_{1k},w-x_{1k}\ra+\la x^*_{1k},x_{1k}-\ox\ra\le\la x^*_{1k},x_{1k}-\ox\ra,\quad k\in\N.
\end{equation}
Note that $\la x^*_{1k},w-\ox\ra\to\la x^*_1,w-\ox\ra$ and $|\la x^*_{1k},x_{1k}-\ox\ra|\le\|x^*_{1k}\|\cdot\|x_{1k}-\ox\|\to 0$ as $k\to\infty$ by the boundedness of $\{x^*_{1k}\}$ in \eqref{ep2}. Passing now to the limit in \eqref{xk} tells us that $\la x^*_1,w-\ox\ra\le 0$ and so $x^*_1\in N(\ox;\O_1)$. It follows from \eqref{ep2} that $-x_1^*\in N(\ox;\O_2)$ and thus $x_1^*\in N(\ox;\O_1)\cap(-N(\ox;\O_2))=\{0\}$, which contradicts the imposed assumption on $\la x^*_{1k}, z\ra\to 0$. Therefore the sequence $\{x^*_{1k}\}$ converges to zero in the weak$^*$ topology of $X^*$, which implies its sequential convergence 
$x^*_{1k}\xrightarrow{w^*}0$ as well. By the assumed SNC property of $\O_1$ at $\ox$ we conclude that $x^*_{1k}\xrightarrow{\|\cdot\|}0$ while yielding $x^*_k \xrightarrow{\|\cdot\|}0$. This surely contradicts \eqref{ep2} and thus ends the proof of implication (ii)$\Longrightarrow$(iii).

To check finally the equivalence assertion in (iii), observe that the separation property \eqref{sep} ensures by Theorem~\ref{extremal principle}(iii) that the sets $\O_1,\O_2$ form an extremal system in $X$, and so we have conditions \eqref{ep1} in (i). Since implication (i)$\Longrightarrow$(ii) has been verified above, this readily justifies the claimed equivalences in (iii) and thus completes the proof of theorem. $\h$

As an immediate consequence of the (convex) approximate extremal principle in Theorem~\ref{convex-ep}(ii), we obtain the celebrated Bishop-Phelps theorem for closed convex sets in general Banach spaces; see \cite[Theorem~3.18]{ph}. Recall that $\ox\in\O$ is a {\em support point} of $\O\subset X$ if there is $0\ne x^*\in X^*$ such that the function $x\mapsto\la x^*,x\ra$ attaints its supremum on $\O$ at $\ox$.\vspace*{-0.1in}

\begin{Corollary}{\bf (Bishop-Phelps theorem).}\label{bp} Let $\O$ be a nonempty, closed, and convex subset of a Banach space $X$. Then the support points of $\O$ are dense on the boundary of $\O$.
\end{Corollary}\vspace*{-0.1in}
{\bf Proof.} It is obvious from \eqref{setex} and the definition of boundary points that for any boundary point $\ox$ of $\O$, the sets $\O_1:=\{\ox\}$ and $\O_2:=\O$ form an extremal system in $X$. Then the result follows from \eqref{ep1} and the normal cone structure in \eqref{nor}.$\h$

Note that a geometric approach involving the approximate extremality conditions \eqref{ep1} at points nearby may be useful for applications to the so-called {\em sequential convex subdifferential calculus} initiated by Attouch-Baillon-Th\'era \cite{abt} and Thibault \cite{thib} in different frameworks and then developed in other publications. Likewise it can be applied as a geometric device of coderivative and conjugate calculus rules, which is our intention in the future research.\vspace*{-0.2in}

\section{Normal Cone Intersection Rule}
\setcounter{equation}{0}\vspace*{-0.1in}

In this section we employ the set extremality and the results of Theorem~\ref{extremal principle} to obtain the exact intersection rule for the normal cone \eqref{nor} under a new qualification condition.

The following theorem justifies a precise representation of the normal cone $N(\ox;\O_1\cap\O_2)$ via normals to each sets $\O_1$ and $\O_2$ under the new qualification condition \eqref{qc} depending on $\ox$, which is weaker than the standard interiority condition in LCTV spaces. For convenience we refer to \eqref{qc} as to the {\em bounded extremality condition}.\vspace*{-0.1in}

\begin{Theorem}{\bf (intersection rule).}\label{nir} Let $\Omega_1,\Omega_2\subset X$ be convex, and let $\ox\in\Omega_1\cap\Omega_2$. Suppose that there exists a bounded convex neighborhood $V$ of $\ox$ such that
\begin{equation}\label{qc}
0\in\mbox{\rm int}\big(\Omega_1-(\Omega_2\cap V)\big).
\end{equation}
Then we have the normal cone intersection rule
\begin{equation}\label{ni}
N(\ox;\Omega_1\cap\Omega_2)=N(\ox;\Omega_1)+N(\ox;\Omega_2).
\end{equation}
\end{Theorem}\vspace*{-0.1in}
{\bf Proof.} To verify \eqref{ni} under the qualification condition \eqref{qc}, denote $A:=\Omega_1$ and $B:=\Omega_2\cap V$ and observe that $0\in\mbox{\rm int}(A-B)$ and $B$ is bounded. Fixing an arbitrary normal $x^*\in N(\ox;A\cap B)$, we get by \eqref{nor} that $\la x^*,x-\ox\ra\le 0$ for all $x\in A\cap B$. Consider the sets
\begin{equation}\label{theta}
\Theta_1:= A\times[0,\infty)\;\;\mbox{\rm and }\;\Theta_2:=\big\{(x,\mu)\in X\times\R\big|\;x\in B,\;\mu\le\la x^*,x-\ox\ra\big\}.
\end{equation}
It follows from the constructions of $\Theta_1$ and $\Theta_2$ that for any $\alpha>0$ we have
\begin{equation*}
\big(\Theta_1+(0,\alpha)\big)\cap\Theta_2=\emp,
\end{equation*}
and thus these sets form an {\em extremal system} by Definition~\ref{ext-sys}. Employing Theorem~\ref{extremal principle}(i) tells us that $0\notin\mbox{\rm int}(\Theta_1-\Theta_2)$. To check next that $\mbox{\rm int}(\Theta_1-\Theta_2)\ne\emp$, take $r>0$ such that $U:=\B(0;r)\subset A-B$. The boundedness of the set $B$ allows us to choose $\bar{\lambda}\in\R$ satisfying
\begin{equation}\label{lambda}
\bar{\lambda}\ge\sup_{x\in B}\la-x^*,x-\ox\ra.
\end{equation}
Then we get ${\rm int}(\Theta_1-\Theta_2)\ne\emp$ by showing that $U\times(\bar{\lambda},\infty)\subset\Theta_1-\Theta_2$. To verify the latter, fix any $(x,\lambda)\in U\times(\bar{\lambda},\infty)$ for which we clearly have $x\in U\subset A-B$ and $\lambda>\bar{\lambda}$, and so $x=w_1-w_2$ with some $w_1\in A$ and $w_2\in B$. This implies in turn the representation
$$
(x,\lambda)=(w_1,\lambda-\bar\lambda)-(w_2,-\bar\lambda).
$$
Further, it follows from $\lambda-\bar\lambda>0$ that $(w_1,\lambda-\bar\lambda)\in\Theta_1$, and we deduce from \eqref{theta} and \eqref{lambda} that $(w_2,-\bar\lambda)\in\Theta_2$, which shows that $\mbox{\rm int}(\Theta_1-\Theta_2)\ne\emp$. Applying now Theorem~\ref{extremal principle}(ii) to the sets $\Theta_1,\Theta_2$ in \eqref{theta} gives us $y^*\in X^*$ and $\gamma\in\R$ such that $(y^*,\gamma)\ne(0,0)$ and
\begin{equation}\label{convexseparation}
\la y^*,x\ra +\lambda_1\gamma\le\la y^*,y\ra+\lambda_2\gamma\;\;\mbox{\rm whenever }\;(x,\lambda_1)\in\Theta_1,\;(y,\lambda_2)\in\Theta_2.
\end{equation}
Using \eqref{convexseparation} with $(\ox,1)\in\Theta_1$ and $(\ox,0)\in\Theta_2$ yields $\gamma\le 0$. Supposing $\gamma=0$, we get
$$
\la y^*,x\ra\le\la y^*,y\ra\;\;\mbox{\rm for all }\;x\in A,\;y\in B.
$$
Since $U\subset A-B$, it readily produces $y^*=0$, a contradiction, which shows that $\gamma<0$. Employing next \eqref{convexseparation} with $(x,0)\in\Theta_1$ for $x\in A$ and $(\ox,0)\in\Theta_2$ tells us that
\begin{equation*}
\la y^*,x\ra\le\la y^*,\ox\ra\;\;\mbox{\rm for all }\;x\in A,\;\;\mbox{\rm and so }\;y^*\in N(\ox;A).
\end{equation*}
Using finally \eqref{convexseparation} with $(\ox,0)\in\Theta_1$ and $(y,\la x^*,y-\ox\ra)\in\Theta_2$ for $y\in B$ implies that
\begin{equation*}
\la y^*,\ox\ra\le\la y^*,y\ra+\gamma\la x^*,y-\ox\ra\;\;\mbox{\rm for all }\;y\in B.
\end{equation*}
Dividing both sides of the obtained inequality by $\gamma<0$, we arrive at
\begin{equation*}
\la x^*+y^*/\gamma,y-\ox\ra\le 0\;\;\mbox{\rm for all }\;y\in B,
\end{equation*}
which verifies by \eqref{nor} the validity of the inclusions
$$
x^*\in-y^*/\gamma+N(\ox;B)\subset N(\ox;A)+N(\ox;B)
$$
and thus shows that $N(\ox;A\cap B)\subset N(\ox;A)+N(\ox;B)$. The opposite inclusion therein is trivial, and so we get the equality $N(\ox;A\cap B)= N(\ox;A)+N(\ox;B)$. Since $N(\ox;A\cap B)=N(\ox;\Omega_1\cap\Omega_2)$ and $N(\ox;B)=N(\ox;\Omega_2)$, it justifies \eqref{ni} and completes the proof.
$\h$\vspace*{-0.1in}

\begin{Remark}{\bf (comparing qualification conditions for the normal intersection formula).}\label{qc-comp} {\rm We have the following useful observations:

{\bf (i)} It is easy to see that, if one of the sets $\O_1,\O_2$ is bounded, the introduced qualification condition \eqref{qc} reduces to the {\em difference interiority condition}
\begin{equation}\label{dqc}
0\in{\rm int}(\O_1-\O_2).
\end{equation}
Furthermore, \eqref{qc} surely holds under the validity of the {\em classical interiority condition} $\O_1\cap({\rm int}\,\O_2)\ne\emp$, which is the only condition previously known to us that ensures the validity of the intersection formula \eqref{ni} in the general LCTV (or even normed) space setting. Indeed, if the latter condition is satisfied, take $u\in\Omega_1\cap({\rm int}\,\Omega_2)$ and $\gg>0$ such that $u+\gg\B\subset\Omega_2$. Then we choose $r>0$ with $u+\gg\B\subset\Omega_2\cap\B(\ox;r)$. Thus $\gg\B\subset\Omega_1-(\Omega_2\cap\B(\ox;r))$ and so $0\in\mbox{\rm int}(\Omega_1-(\Omega_2\cap V))$, where $V:=\B(\ox;r)$.

As the following simple example shows, the bounded extremality condition \eqref{qc} may be weaker than the classical interiority condition even in $\R^2$. Indeed, consider the convex sets
$$
\O_1:=\R\times[0,\infty)\;\;\mbox{\rm and }\;\O_2:=\{0\}\times\R
$$
for which $\O_1\cap({\rm int}\,\O_2)=\emp$, while the conditions $0\in{\rm int}(\O_1-\O_2)$ and \eqref{qc} hold.

{\bf (ii)} If $X$ is {\em Banach} and both sets $\O_1,\O_2$ are {\em closed} with $\mbox{\rm int}(\O_1-\O_2)\ne\emp$, the difference interiority condition \eqref{dqc} reduces to Rockafellar's {\em core qualification condition} $0\in{\rm core}(\O_1-\O_2)$ introduced in \cite{r1}. This follows from the equivalence
\begin{equation}\label{core-int}
\big[0\in{\rm core}(\O_1-\O_2)\big]\Longleftrightarrow\big[0\in{\rm int}(\O_1-\O_2)\big]
\end{equation}
valid in this case. Indeed, the implication ``$\Longleftarrow$" in \eqref{core-int} is obvious due to $\i\O\subset\core\O$ for any set. To verify the opposite implication in \eqref{core-int}, recall that $\i\O=\core\O$ for closed convex subsets of Banach spaces by \cite[Theorem~4.1.8]{BZ}. Using now the well-known fact that $\i\Bar\O=\i\O$ for convex sets with nonempty interiors yields
$$
0\in\mbox{\rm core}\big(\Bar{\O_1-\O_2}\big)=\mbox{\rm int}\big(\Bar{\O_1-\O_2}\big)=\mbox{\rm int}\big(\O_1-\O_2\big).
$$
Note that the core qualification condition is superseded  in the same setting by the requirement that $\R^+(\Omega_1-\Omega_2)\subset X$ is a closed subspace, which is known as the {\em Attouch-Br\'ezis regularity condition} established in \cite{AB} with the usage of convex duality and the fundamental Banach-Dieudonn\'e-Krein-\v Smulian theorem in general Banach spaces.}
\end{Remark}\vspace*{-0.1in}

The next proposition shows that the core condition $0\in{\rm core}(\O_1-\O_2)$ implies the extremality one \eqref{qc} for closed subsets of reflexive Banach spaces {\em provided that} ${\rm int}(\O_1-\O_2)\ne\emp$. Thus the extremality approach of Theorem~\ref{nir} offers in this setting a simplified proof of the intersection formula in comparison with those known in the literature.\vspace*{-0.1in}

\begin{Proposition}{\bf (bounded extremality condition in reflexive spaces).}\label{intersection rule reflexive} The qualification condition \eqref{qc} holds at any $\ox\in\Omega_1\cap\Omega_2$ if $X$ is a reflexive Banach space and $\Omega_1,\Omega_2\subset X$ are closed convex sets such that ${\rm int}(\O_1-\O_2)\ne\emp$ and $0\in\mbox{\rm core}(\Omega_1-\Omega_2)$.
\end{Proposition}\vspace*{-0.1in}
{\bf Proof.} Fix any number $r>0$ and show that
\begin{equation}\label{cl}
0\in\mbox{\rm core}\big(\Omega_1\cap\B(\ox;r)-\Omega_2\cap\B(\ox;r)\big).
\end{equation}
Indeed, the assumption ${\rm int}(\O_1-\O_2)\ne\emp$ allows us to find $\gamma>0$ such that $\gamma\B\subset\Omega_1-\Omega_2$. For any $x\in X$ denote $u:= \frac{\gamma}{\|x\|+1}x\in\gamma\B$ and get $u=w_1-w_2$ with $w_i\in\Omega_i$ for $i=1,2$. Hence there is a constant $\bar{\gamma}>0$ depending on $x$ and $r$ for which
$$
t\max\big\{\|w_1-\ox\|,\|w_2-\ox\|\big\}<r\;\mbox{ whenever }\;0<t<\bar{\gamma}.
$$
This readily justifies the relationships
\begin{equation*}
tu=tw_1-tw_2=\big(\ox+t(w_1-\ox)\big)-\big(\ox+t(w_2-\ox)\big)\in\big(\Omega_1\cap\B(\ox;r)\big)-\big(\Omega_2\cap\B(\ox;r)\big)
\end{equation*}
for all $0<t<\bar{\gamma}$ and thus establishes the claimed inclusion \eqref{cl} by the core definition \eqref{core-def}.

Since $X$ is reflexive and the sets $\Omega_i\cap\B(\ox;r)$, $i=1,2$, are closed and bounded in $X$, they are weakly sequentially compact in this space. This implies that their difference $\big(\Omega_1\cap\B(\ox;r)\big)-\big(\Omega_2\cap\B(\ox;r)\big)$ is closed in $X$. Then we get by \cite[Theorem~4.1.8]{BZ} that
$$
0\in\mbox{\rm core}\big(\Omega_1\cap\B(\ox;r)-\Omega_2\cap\B(\ox;r)\big)=\mbox{\rm int}\big(\Omega_1\cap\B(\ox;r)-\Omega_2\cap\B(\ox;r)\big)\subset{\rm int}\big(\Omega_1-\Omega_2\cap\B(\ox;r)\big),
$$
which verifies \eqref{qc} and thus completes the proof of the proposition. $\h$\vspace*{-0.2in}

\section{Support Functions for Set Intersections}
\setcounter{equation}{0}\vspace*{-0.1in}

In this section we derive a precise representation of support functions for convex set intersections via the infimal convolution of the support functions to the intersection components under the {\em difference interiority condition} \eqref{dqc}. This result under \eqref{dqc} seems to be new in the literature on convex analysis in LCTV (and also in normed) spaces; see Remark~\ref{rem-supp} for more discussions. Furthermore, we present a novel geometric device for results of this type by employing set extremality and the normal intersection rule obtained above.

Recall that the {\em support function} of a nonempty set $\O\subset X$ is given by
\begin{equation}\label{supp}
\sigma_\O(x)(x^*):=\sup\{\la x^*,x\ra\big|\;x\in\O\big\},\quad x^*\in X^*.
\end{equation}
The {\em infimal convolution} of two functions $f,g\colon X\to\oR$ is
\begin{equation}\label{ic}
(f\oplus g)(x):=\inf\big\{f(x_1)+g(x_2)\big|\;x_1+x_2=x\big\}=\inf\big\{f(u)+g(x-u)\big|\;u\in X\big\}.
\end{equation}\vspace*{-0.35in}

\begin{Theorem}{\bf(support functions for set intersections via infimal convolutions).}\label{sigma intersection rule} Let the sets $\Omega_1,\Omega_2\subset X$ be nonempty and convex, and let and one of them be bounded. Then the difference interiority condition \eqref{dqc} ensures the representation
\begin{equation}\label{convol}
(\sigma_{\Omega_1\cap\Omega_2})(x^*)=(\sigma_{\Omega_1}\oplus\sigma_{\Omega_2})(x^*)\;\;\mbox{\rm for all }\;x^*\in X^*.
\end{equation}
Moreover, for any $x^*\in\dom(\sigma_{\Omega_1\cap\Omega_2})$ there are $x^*_1,x^*_2\in X^*$ such that $x^*=x^*_1+x^*_2$ and
\begin{equation}\label{convol1}
(\sigma_{\Omega_1\cap\Omega_2})(x^*)=\sigma_{\Omega_1}(x^*_1)+\sigma_{\Omega_2}(x^*_2).
\end{equation}
\end{Theorem}\vspace*{-0.1in}
{\bf Proof.} First we check that the inequality ``$\le$" in \eqref{convol} holds in the general setting. Fix any $x^*\in X^*$ and pick $x^*_1,x^*_2\in X^*$ such that $x^*=x^*_1+x^*_2$. Then it follows from \eqref{supp} that
\begin{equation*}
\la x^*,x\ra=\la x^*_1,x\ra +\la x^*_2,x\ra\le\sigma_{\Omega_1}(x^*_1)+\sigma_{\Omega_2}(x^*_2)\;\;\mbox{\rm whenever }\;x\in\Omega_1\cap\Omega_2.
\end{equation*}
Taking the infimum on the right-hand side above with respect to all $x^*_1,x^*_2\in X^*$ satisfying $x^*_1+x^*_2=x^*$ gives us by definition \eqref{ic} of the infimal convolution that
\begin{equation*}
\la x^*,x\ra\le(\sigma_{\Omega_1}\oplus\sigma_{\Omega_2})(x^*).
\end{equation*}
This verifies and the inequality ``$\le$" in \eqref{convol} by taking the supremum on the left-hand side therein with respect to $x\in\Omega_1\cap\Omega_2$.

To justify further the opposite inequality in \eqref{convol} under the validity of \eqref{dqc}, suppose that $\Omega_2$ is bounded. It suffices to consider the case where $x^*\in\dom(\sigma_{\Omega_1\cap\Omega_2})$ and prove the inequality ``$\le$" in \eqref{convol1}; then the one in \eqref{convol} and both statements of the theorem follow.

To proceed, denote $\alpha:=(\sigma_{\Omega_1\cap\Omega_2})(x^*)\in\R$, for which we clearly have $\langle x^*,x\rangle-\alpha\le 0$ whenever $x\in\Omega_1\cap\Omega_2$, and then construct the two nonempty convex subsets of $X\times\R$ by
\begin{equation}\label{Theta}
\Theta_1:=\Omega_1\times[0,\infty)\;\;\mbox{\rm and }\;\Theta_2:=\big\{(x,\lambda)\in X\times\R\big|\;x\in\Omega_2,\;\lambda\le\langle x^*,x\rangle-\alpha\big\}.
\end{equation}
Observe that the sets $\Th_1,\Th_2$ form an {\em extremal system}. Indeed, it follows from the choice of $\al$ and the construction in \eqref{Theta} that
for any $\gamma>0$ we have
\begin{equation*}
\big(\Theta_1+(0,\gamma)\big)\cap\Theta_2=\emp.
\end{equation*}
Then Theorem~\ref{extremal principle}(i) tells us that $0\notin\mbox{\rm int}(\Theta_1-\Theta_2)$. Arguing similarly to the proof of Theorem~\ref{nir}, we see that the condition $\mbox{\rm int}(\Theta_1-\Theta_2)\ne\emp$ holds for the sets in \eqref{Theta}. Thus Theorem~\ref{extremal principle}(ii) allows us to find a pair
$(y^*,\beta)\ne(0,0)$ such that
\begin{equation}\label{sep-conv}
\la y^*, x\ra+\lambda_1\beta\le\la y^*,y\ra+\lambda_2\beta\;\;\mbox{\rm whenever }\;(x,\lambda_1)\in\Theta_1,\;(y,\lambda_2)\in\Theta_2.
\end{equation}
Choosing $(\ox,1)\in\Theta_1$ and $(\ox,0)\in\Theta_2$ in \eqref{sep-conv} shows that $\beta\le 0$. If $\beta=0$, then
\begin{equation*}
\la y^*,x\ra\le\la y^*,y\ra\;\;\mbox{\rm for all }\;x\in\Omega_1,\;y\in\Omega_2.
\end{equation*}
By ${\rm int}(\Omega_1-\Omega_2)\ne\emp$ this yields $y^*=0$, a contradiction justifying the negativity of $\bb$ in \eqref{sep-conv}. Take now $(x,0)\in\Theta_1$ and $(y,\langle x^*,y\rangle-\alpha)\in\Theta_2$ in \eqref{sep-conv} and then get
\begin{equation*}
\la y^*,x\ra\le\la y^*,y\ra+\beta(\la x^*,y\ra-\alpha),
\end{equation*}
which can be equivalently rewritten (due to $\bb<0$) as
\begin{equation*}
\alpha\ge\big\la y^*/\beta+x^*,y\big\ra+\big\la-y^*/\bb,x\big\ra\;\;\mbox{\rm for all }\;x\in\Omega_1,\;y\in\Omega_2.
\end{equation*}
Denoting $x^*_1:=y^*/\beta+x^*$ and $x^*_2:=-y^*/\beta$, we have $x^*_1+x^*_2=x^*$ and $\langle x^*_1,x\rangle +\langle x^*_2,y\rangle\le\alpha$ for all
$x\in\Omega_1$ and $y\in\Omega_2$. This shows that
\begin{equation*}
\sigma_{\Omega_1}(x^*_1)+\sigma_{\Omega_2}(x^*_2)\le\alpha=\sigma_{\Omega_1\cap\Omega_2}(x^*)
\end{equation*}
and thus completes the proof of the theorem.$\h$\vspace*{-0.1in}

\begin{Remark}{\bf(comparison with Fenchel duality).}\label{rem-supp}
{\rm Since the qualification condition \eqref{qc} used in Theorem~\ref{nir} is equivalent to \eqref{dqc} employed in Theorem~\ref{sigma intersection rule} when one of the sets $\O_1,\O_2$ is bounded, all the comments given in Remark~\ref{qc-comp} are applied here. On the other hand, there is a remarkable feature of the calculus rules for support functions presented in Theorem~\ref{sigma intersection rule}, which does not have analogs in the setting of Theorem~\ref{nir} and should be specially commented. Namely, the support function \eqref{supp} is the {\em Fenchel conjugate}
\begin{equation*}
f^*(x^*):=\sup\big\{\la x^*,x\ra-f(x)\big|\;x\in X\big\},\quad x^*\in X^*,
\end{equation*}
of the indicator function $f(x):=\dd(x;\O)$ of a given set $\O\subset X$, and hence a well-developed {\em conjugate calculus} can be applied to establish representations \eqref{convol} and \eqref{convol1}; see, e.g., the books \cite{Ra,r1,s,z} and the references therein. However, it seems to us that such an approach from Fenchel duality misses the specific results of Theorem~\ref{sigma intersection rule} derived for the support function under the qualification condition \eqref{dqc} in general LCTV spaces. Observe also that, in contrast to analytical schemes usually applied to deriving conjugate calculus and then deducing results of the type of Theorems~\ref{nir} and \ref{sigma intersection rule} from them, we develop here a {\em geometric approach} in the other direction based on {\em set extremality}.}
\end{Remark}


\small
\begin{thebibliography}{99}

\bibitem{abt} Attouch H., Baillon, J.-B., Th\'era, M.: Variational sum of monotone operators. J. Convex Anal. 1, 1--29 (1994)

\bibitem{AB} Attouch, H., Br\'ezis, H.: Duality of the sum of convex functions in general Banach spaces. In J. A. Barroso, J.A. (ed.) Aspects of Mathematics
and Its Applications 34, pp.\ 125–-133. North-Holland, Amsterdam (1986)

\bibitem{bc} Bauschke, H.H., Combettes, P.L.: Convex Analysis and Monotone Operator Theory in Hilbert Spaces. Springer, New York (2011)

\bibitem{bl} Borwein, J.M., Lewis, A.S.: Convex Analysis and Nonlinear Optimization, 2nd edition. Springer, New York (2006)

\bibitem{blm} Borwein, J.M., Lucet, Y., Mordukhovich, B.S.: Compactly epi-Lipschitzian convex sets and functions in normed spaces. J. Convex Anal. 7,
375--393 (2000)

\bibitem{BS} Borwein, J.M., Str\'ojwas H.: Tangential approximations. Nonlinear Anal. 9, 1347--1366 (1985)

\bibitem{BZ} Borwein, J.M., and Zhu, Q.J.: Techniques of Variational Analysis. Springer, New York (2005)

\bibitem{Ra} Bo\c t, R.I.: Conjugate Duality in Convex Optimization. Springer, Berlin (2010)

\bibitem{c} Clarke, F.H.: Optimization and Nonsmooth Analysis. Wiley, New York (1983)

\bibitem{fm} Fabian, M., Mordukhovich, B.S.: Sequential normal compactness versus topological normal compactness in variational analysis.
Nonlinear Anal. 54, 1057--1067 (2003)

\bibitem{HU} Hiriart-Urruty, J.-B., Lemar\'echal, C.: Convex Analysis and Minimization Algorithms I, II. Springer, Berlin (1993)

\bibitem{m-book1} Mordukhovich, B.S.: Variational Analysis and Generalized Differentiation, I: Basic Theory, II: Applications. Springer, Berlin (2006)

\bibitem{bmn} Mordukhovich, B.S., Nam, N.M.: An Easy Path to Convex Analysis and Applications. Morgan \& Claypool Publisher, San Rafael, CA (2014)

\bibitem{bmn1} Mordukhovich, B.S., Nam, N.M.: Geometric approach to convex subdifferential calculus. Optimizaion (2016). doi: 10.1080/02331934.2015.1105225

\bibitem{ph} Phelps, R.R.: Convex Functions, Monotone Operators and Differentiability, 2nd edition. Springer, Berlin (1993)

\bibitem{r} Rockafellar, R.T.: Convex Analysis. Princeton University Press. Princeton, NJ (1970)

\bibitem{r1} Rockafellar, R.T.: Conjugate Duality and Optimization. SIAM, Philadelphia, PA (1974)

\bibitem{RockWets-VA} Rockafellar, R.T., Wets, R.J-B.: Variational Analysis. Springer, Berlin (1998)

\bibitem{s} Simons, S.: From Hahn-Banach to Monotonicity, 2nd edition. Springer, Berlin (2008)

\bibitem{thib} Thibault, L.: Sequential convex subdifferential calculus and sequential Lagrange multipliers. SIAM J. Control Optim. 35,
1434--1444 (1997)

\bibitem{z} Z\u{a}linescu, C.: Convex Analysis in General Vector Spaces. World Scientific, Singapore (2002)
\end{thebibliography}
\end{document}